%% file: main.tex
\documentclass[runningheads]{llncs}
\usepackage[T1]{fontenc}
\usepackage{graphicx}
\usepackage[hidelinks]{hyperref}
\usepackage[justification=centering]{caption}
\input{packages.tex}

\input{DeclareMathOperator.tex}

\input{preamble_ru.tex}

\input{newcommands.tex}

\newcommand{\bx}{\bar x}
\newcommand{\bA}{\bar A}
\newcommand{\GLPK}{{\tt GLPK}\xspace}
\newcommand{\MATLAB}{{\tt MATLAB}\xspace}
\newcommand{\OCTAVE}{{\tt OCTAVE}\xspace}

\begin{document}
\title{Linear Optimization by Conical Projection\thanks{Research was supported by the subsidy of the Ministry of Science and High Education of Russian Federation 075-02-2023-946 from February
16, 2023.}}
%
%
\author{Evgeni Nurminski\inst{1},
Roman Tarasov\inst{2} }
\authorrunning{Е. Nurminski et al.}
%
\institute{Far Eastern Centre for Research and Education in Mathematics, Far Eastern Federal University, Russia \\
\email{nurminskiy.ea@dvfu.ru}\\
 \and
Skolkovo Institute of Science and Technology, Russia\\
\email{troman85235@gmail.com}}
\maketitle              

\renewcommand{\abstractname}{}
\begin{abstract}
This article focuses on numerical efficiency of projection algorithms for solving linear optimization problems.
The theoretical foundation for this approach is provided by the basic result that bounded finite dimensional linear optimization problem can be solved by single projection operation on the feasible polyhedron. The further simplification transforms this problem into projection of a special point onto a convex polyhedral cone generated basically by inequalities of the original linear optimization problem.  


\end{abstract}
\section*{Introduction}
Linear optimization remains a workhorse of many practical applications and
modern general-purpose industrial quality simplex-based algorithms and interior point methods
demonstrated remarkable success in the area. 
Nevertheless there is an intellectual challenge to develop new approaches which may find their applications
in these or that situations. 
In this article, we consider a projection-based algorithm for solving linear optimization problems in the standard form
\bleq{stanlp}{
\min_{Ax \leq b} cx ~ = ~c \xs = ~\max_{u A = c,~ u \leq 0} ub ~ = ~ u^\star b
}
which is written here with its dual.

Here in more or less standard notation \(x\), \(c\) and \(\xs\) are vectors of finite (\(n\)) dimensional euclidean space \(E\) with the inner product \(xy\) and norm \(\|x\|^2 = x x\). The right-hand side \(b\)
of the constraints in (\ref{stanlp}) belongs to \(m\)-dimensional space \(E'\) and \(A\) is a linear operator (\(m\times n\) matrix) from \(E\) to \(E'\) --- \(A: E \to E'\).
The dual part of this problem has \(m\)-dimentional vector \(u\) of dual variables
and its optimal solution is denoted as \(u^\star\).

Solution of (\ref{stanlp}) can be reduced to solving primal-dual system of linear inequalities
or convex and even polyhedral feasibility problem (PFP):
\bleq{prim-dual-system}{ cx \geq bu, ~ Ax \leq b, ~ uA = c, ~u \leq 0.}
The PFP and Convex Feasibility Problem (CFP) for the general convex sets were
in gun sights of many mathematicians since the middle of 20 century
and projection methods are amongst the most popular for solving it since pioneering works
\cite{kacz37},
\cite{Cimmino38},
\cite{neuman50}
(see the extensive review 
of H. Bauschke and J.M. Borwein \cite{babo96}).
However in the area of linear optimization projection methods were not very successful 
in practical sense, mainly for slow  convergence and  computational difficulties of solving multiple high dimensional projection problems for polyhedrons of the general type \cite{nogood}.

We hope nevertheless that our experiments inspire new interest in projection agenda and algorithms. These hopes are based on certain intrinsic properties of projection operation which can operate on reduced basis sets and on additional decomposition possibilities, see f.i. \cite{good}. 

\makeatletter{}\section{Notations and Preliminaries}
\label{prenot}
As it defined in Introduction
let \(E\) be a finite-dimensional vector space of the primal variables
with the standard inner product \(xy\) and the norm \(\|x\|^2 = xx\).
This space is then self-conjugate with the duality relation induced by the inner product.
The dimensionality of this space, if needed, is determined as \(\dim(E)\) and the space of dimensionality \(n\)
when necessary is denoted as \(E^n\). 
The non-negative part of a space \(E\) will be denoted as \(E_+\).

Among the others special vectors and sets we mention the null vector \(\vnll\),
vector of ones \(\vone = (1,1, \dots, 1)\), and the standard simplex \(\Delta = E_+ \cap \lbrace x: \vone x = 1 \rbrace \). 
Linear envelope, convex and conical hull of a set \(X\) are denoted as \(\lin(X)\), \(\co(X)\) and \(\Co(X)\)
respectively.

We define linear operators, acting from \(E\) into \(E'\) with \(\dim(E') = m \)
as collections of vectors \( {\cal A} = \lbrace a^1, a^2, \dots, a^m \rbrace ~ \mbox{with} ~ a^i \in E \) which produce vector
\(y = (y_1, y_2, \dots, y_m) \in E^m \) according to
following relations \(y_i = a^i x, i = 1, 2, \dots, m \).
In the classical matrix-vector notation vectors \( {\cal A} \) form the rows of the matrix \(A\) and \(y = A x \).
At the same time we will consider the row subspace \(E'\) as the linear envelope
of \(\cal A\):
\[ E' = \lin ({\cal A}) = \{ x = \sum_{i = 1}^m a^i z_i = A^T z, z \in E^m  \} \subset E. \]

The projection operator of a point \(p\) onto a closed convex set \(X\) in \(E\) is defined as
$$ \prc{p}{X} = \amin_{x \in X} \| p - x \|, $$
that is \( \min_{x \in X} \| p - x \| = \| p - \prc{p}{X} \|\).
For closed convex \(X\), this operator is well-defined and Lipschitz-continuous with the Lipschitz constant less or equal \(1\).
The point-to-set projection operation is naturally generalized for sets:
\( \prc{X}{A} = \{ z = \prc{x}{A}, x \in X \}\). 

We will also notice that this operator is idempotent:
\(\prc{(\prc{p}{X})}{X} = \prc{p}{X}\) and linear for projection on
linear subspace \(L\) of \(E\):
\( \prc{\alpha p }{L} = \alpha (\prc{p}{L}) \) for \(\alpha \in \R \) and \(\prc{(p + q)}{L} = \prc{p}{L} + \prc{q}{L} \).
Of course \(p = \prc{p}{L} + \prc{p}{L^\perp}\).

For a closed convex set \(X\) denote as \(\supp{X}{z}\) its support function
\begin{equation}\label{supf}
 (X)_z = \begin{array}[t]{c} \min \\ x \in X \end{array} xz. 
\end{equation}
In this notation the standard linear optimization problem
\begin{equation}\label{lopt}
\minar{Ax \leq b} ~cx ~=~ \minar { x \in X } cx
\end{equation}
becomes just \( (X)_c\).

Basically the same holds and for nonlinear problems
\begin{equation}\label{nlopt}
\minar{x \in X} ~ f(x) ~=~ \minar { \bx \in \bar X } \bar c \bx = \supp{\bar X}{\bar c}
\end{equation}
for \(\bar x = (x, \xi) \), \(\bar X = \{ \bx : x \in X, f(x) \leq \xi\), \(\bar c = (\vnll, 1) \).

There is a general result which connect support functions with projection \cite{asiq}.
\begin{theorem}
\label{thm1}
Let \(X\) --- closed bounded subset of \(E\) and \(c \in E \).
Then for any \(x^0\)
\begin{equation}
\label{supfun}
\supp{X}{c} = \lim_{\tau \to \infty} ~ c (\prc{(x^0 + \tau c)}{X}). 
\end{equation}
\end{theorem}
For the formal correctness of application of the theorem \ref{thm1} to the set \(\bar X\) it is necessary to ensure
boundness of \(\bar X\).
This, generally speaking, formal requirement can be easily satisfied
by adding an arbitrary upper bound
\(\bar f \geq \inf\limits_{x \in X} f(x)\) for \(\xi\).
Toward this purpose any \(x^0 \in X\) will provide
trivial upper bound \(\bar f = f(x^0)\).



For linear optimization problems
(\ref{lopt})
where \(X\) is a bounded polyhedron,
exact equivalence can be proved \cite{fin-projection}:
\begin{theorem} 
If (\ref{lopt}) has a unique solution \(\xs\),
then for any \(x^0\) there existsт \(\theta_c > 0\) such that
\bleq{prl}{\prc{(x^0 - \theta c)}{X} = \xs}
for any \(\theta \geq \theta_c\).
\end{theorem}

In more details the problem (\ref{prl}) can be written down as:
\bleq{x}{
\begin{array}[t]{c} \min \\ x \in X \end{array} \| x - x^0 + \theta c \|^2 =
\begin{array}[t]{c} \min \\ y \in X_{\theta, c} \end{array} \| y \|^2,
}
where
\( X_{\theta, c} = X - x^0 + \theta c \)
is the original feasible set \(X\), shifted by \(x^c = \theta c - x^0\).

If the polyhedron \(X\) is described by a system of linear inequalities
\bleq{lineq}{ X = \{ x: Ax \leq b \} }
then
\bleq{nlineq}{ X_{\theta, c} = \{ x: Ax \leq b^c \}, }
\( b^c = b - A(x^0 - \theta c). \)

The latter problem (\ref{x}) does not look as  something essentially different however
it can be transformed into the conical projection problem,
\bleq{conic}{
\begin{array}[t]{c} \min \\ \bar x \in \Cone (\bar A) \end{array} \| \bar x - p \|^2
}
where \(\bx = (x, \xi)\) --- is the vector from extended space \(\bar E = E \times \R \),
and \( \bA = \vert A, -b^c \vert \).
The rows of this matrix can be considered as vector of
\(\bar E\).
Then \(\Cone(\bar A)\) is the conical envelope of these vectors which can
be represented as
\[ \Cone(\bar A) = \{ \bar x = \bar A^T z, z \in E^r_+ \}, \]
where \(E^r_+\) --- non-negative orthant of the correspondent dimensionality.

Finally the algorithm for solution of (\ref{lopt}) can be represented by
the algorithmic scheme \ref{proalgo}.
\begin{algorithm}
\caption{Solving (\ref{lopt}) by projection.}
\SetAlgoLined
\KwData{ The dataset \((A, b, c\) of the original problem, and scaling constant \( \theta > 0 \)}
\KwResult{The solution \(\xs\) of the linear optimization problem (\ref{lopt}).}
  
{\bf Step 1.} Data preparation for projection problem\;
\bleq{prep}{
\begin{array}{c}
x^c = x - 0c; ~ b^c = b - A x^c;
\\
\bar A = \left[ A, -b^c \right]; ~ \bar p = [\null_n, 1 ]
\end{array}
}
{\bf Step 2.} Solution of the projection problem\;
\bleq{propro}{
\begin{array}[t]{c} \min \\ \bar x \in \Cone(\bar A) \end{array} \| \bar x - \bar p \|^2 =
\| \prc{\bar p}{\Cone(\bar A)} - \bar p \|^2
}
{\bf Step 3.} Getting back to (\ref{lopt})\;
By representing solution of the problem (\ref{propro}) as
\(\prc{\bar p}{\Cone(\bar A)} = (y^c, \xi)\),
where \(y^c \in E\), а \(\xi \in \R\)
compute
\[ \xs = y^c/\xi + \theta c. \] 
\label{proalgo}
\end{algorithm}
\section{Inside out}
The subject of this section is the least norm problem
\(\min_{x \in X} \|x\|^2\)
in an \(n\)-dimensional euclidean space \(E\)
for a bounded closed convex polyhedron \(X\).
Here we do not make a great distinction between row and column vectors which are assumed of
any type depending on context.
Polyhedron \(X\) most commonly
described as the intersection of half-spaces
\bleq{hs}{ X = \{ x: a^i x \leq \beta_i, i = 1,2, \dots, m \} = \{x : Ax \leq b\} }
where vectors \(a^i, i = 1,2,\dots, m\) of the dimensionality \(n\) can be considered as rows of the matrix \(A\),
and the \(m\)-vector \(b = (\beta_1, \beta_2, \dots, \beta_m) \) is the corresponding right-hand side vector.
It can be considered as the ''outer'' description of \(X\) in contrast with the ''inner'' description
\bleq{coex}{ X = \co( \hat x^j, \hat x^j \in \mathrm{Ext} (X), j = 1,2, \dots, J)}
as the convex hull of the set \(\mathrm{Ext}(X)\) of extreme points of the same set \(X\).
The later is often considered as ''polytope'' description.
These are equivalent descriptions for this class of polyhedrons/polytopes,
but direct conversion between them is complicated as any of them may be exponentially long
even for the polynomially long in \(n, m\) counterparts.

The polyhedron description is more common so the vast majority of computational algorithms is developed
namely for this description of \(X\).
The notable exceptions are possibly game problems with probability simplexes and nondifferentiable optimization
algorithms in which subdifferentials are approximated by convex hulls of known subgradients.
However convex hull-like description has its own computational advantages, for instance as linear optimization problem over
convex hulls it has low \(nm\) complexity for the trivial direct algorithm and can be reduced to logarithmic complexity
if parallel computations allowed.
In
\cite{outer_projection}
we considered the transformation of the least norm problem with
the polyhedral description (\ref{hs}) into
the close relative of (\ref{coex}) with practically the same data-size as (\ref{hs}). 
The original version of this transformation was rather convoluted
and here we present its alternative derivation which
uses basically only standard duality arguments.
 
To begin with we expand our basic space \(E\) with one additional variable into \(\bar E = E \times \R\) and
transform the initial least norm problem into something which is almost homogeneous:
\bleq{lis}{
\minar{ A x \leq b } \half (\| x \|^2 + 1) =
\minar{ \bar A \bx \leq 0 \\ \bar e \bx = 1 } \half \| \bx \|^2 }
with \( \bar A = \vert A, -b \vert \), \(\bar x = (x, \xi) \), and vector \(\bar e = (0, 0, \dots, 0, 1) \in \bar E\). 
The saddle point reformulation of this problem goes as follows:
\[
\begin{array}{c}
\minar{\bA \bx \leq 0 \\ \bar e \bx = 1 } \half \| \bx \|^2 = \half \| \bx^\star \|^2 =
\maxar{u \geq 0, ~ \theta}
\minar{\bx } \{\half \| \bx \|^2 + \theta (1 - \bar e \bx) + u \bA \bx \} =
\\
\maxar{u \geq 0, ~ \theta } \{ \theta + 
\minar{\bx } \{ \half \| \bx \|^2 + ( u \bA - \theta \bar e) \bx \} =
\maxar{ u \geq 0, ~ \theta} \{ \theta - \half \| u \bA - \theta \bar e \|^2 \}
\end{array}
\]
Introducing the cone \(\cK = \{ z: z = u \bA, u \geq 0 \}\) we can rewrite the last problem as
\[
-\minar{\theta} \{ \half \min_{z \in \cK} \| z - \theta \bar e \|^2  - \theta\} =
-\minar{\theta} \{ \half \theta^2 \min_{z \in \cK} \| z - \bar e \|^2 - \theta \} =
-\minar{\theta} \{ \half \gamma^2 \theta^2  - \theta\} = 1/2\gamma^2,
\]
where we made use of \( \alpha \cK = \cK \) for any \(\alpha > 0\) and denoted \(\gamma^2 = \min_{z \in \cK} \| z - \bar e \|^2 \).
The solution of the last minimum is attained for \(\theta^\star = 1/\gamma^2\).
As solution of \(\min_{z \in \cK} \| z - \bar e \|^2 = \| z^\star - \bar e \|^2 \) is unique we obtain
\(\bx^\star = \theta^\star (z^\star - \bar e ) \). 
\section{Numerical experiments}
For numerical experiments we implemented the projection algorithm in the \OCTAVE open-source MATLAB-like system
\cite{oct}.
For comparison, we used the \GLPK linear programming kit
implemented in C programming language and also built into \OCTAVE as internal function.
This function however has a limited functionality in comparison with stand-along \GLPK solver
\cite{glpk},
but still allows for some basic comparison of projection algorithm and contemporary symplex method.

We considered the following 2 types of linear optimization problems.
The first one consists of linear optimization problems of different dimensions with the following structure
\bleq{test-lp}{
\minar{Ax \leq \vnll \\ -g \leq x \leq f } c x ~ = \minar {\bA x \leq \bar b } c x
}
where \(\bA = [ A; ~ I; ~ -I ], ~ \bar b = [ \vnll; ~ f; ~ g ] \).
Here \(I\) is the identity matrix,
the elements of the matrix \(A\) and the vectors \( f, g \)
are randomly generated independently and uniformly from the segment \([0,1]\),
and the elements of the vector \(c\) are generated from the segment \([-5 ,5]\).
Also for matrices  \( A, B ) \) of matching dimensions \( [ A; B ] \) denotes
(by following \MATLAB/\OCTAVE convention) a stacked up \(A, B\),
that is \( [ A^T B^T ]^T\).
This form of problems on one hand ensures the feasibility ( \(\vnll\) is a feasible point ) and
boundness (due to simple two-sided bounds on variables) but  intersection of the cone and the box produces sufficiently demanding feasible polyhedron.

However, this form provides one-sided advantages to \GLPK as it immediately makes use of built-in presolver which takes into account sparsity and specific structure of box constraints.
Of course our experimental implementation of projection algorithm is missing such advanced features,
so our second type of test problems was a simple modification of (\ref{test-lp}) 
which consisted in
replacing \(x\) with the new variable \(z\), such that \(x = Qz\)
with 100\%-dense random unitary matrix \(Q\).
Then the test problems become
\bleq{test-lp-Q}{
\minar {\bA_Q z \leq \bar b } c_Q z
}
where \(\bA_Q = \bA Q, ~ c_Q = c Q.\)
After such changes the problem constraints became fully dense and \GLPK presolver
does not interfere with optimization. 

\begin{figure}[h]
\includegraphics{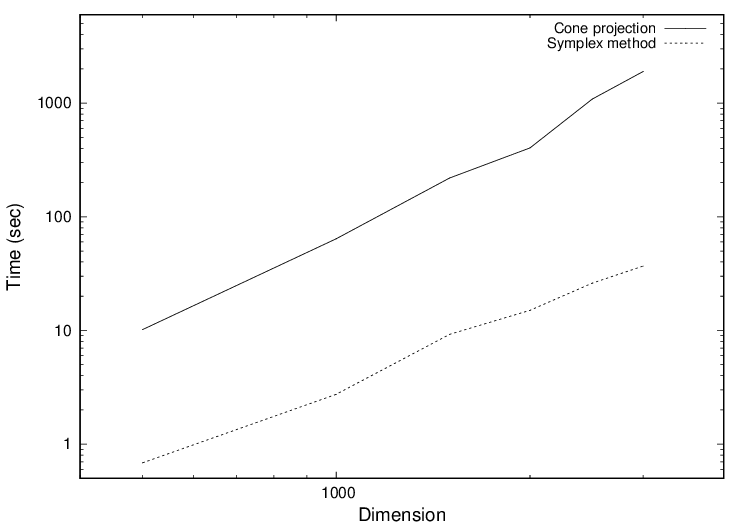}
\caption{Dependence of the running time of the algorithm on the dimensionality: \eqref{test-lp}. }
\label{test_lp1}
\end{figure}

Fig.\ref{test_lp1} shows that presolver keeps the \GLPK solver noticeably faster, but projection algorithm demonstrates at least the similar dynamics when problem dimensions increase.

\begin{figure}[h]
\includegraphics{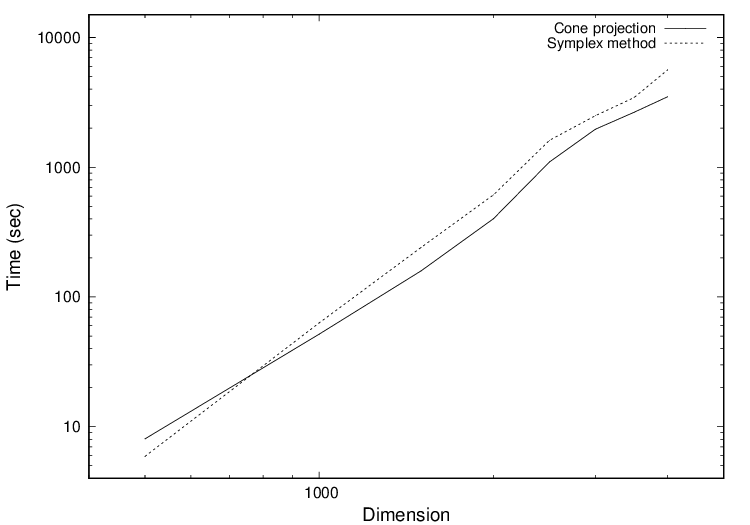}
\caption{Dependence of the running time of the algorithm on the dimensionality: \eqref{test-lp-Q}.}
\label{test_lp2}
\end{figure}



For dense problems \eqref{test-lp-Q} projection algorithm demonstrates however some speed-up
(Fig. \ref{test_lp2})
which slightly increases when problem size grows.
It was a pleasant surprise that despite the very different levels of implementation 
the projection algorithm was faster than \GLPK.

\begin{figure}
\includegraphics[width=\textwidth]{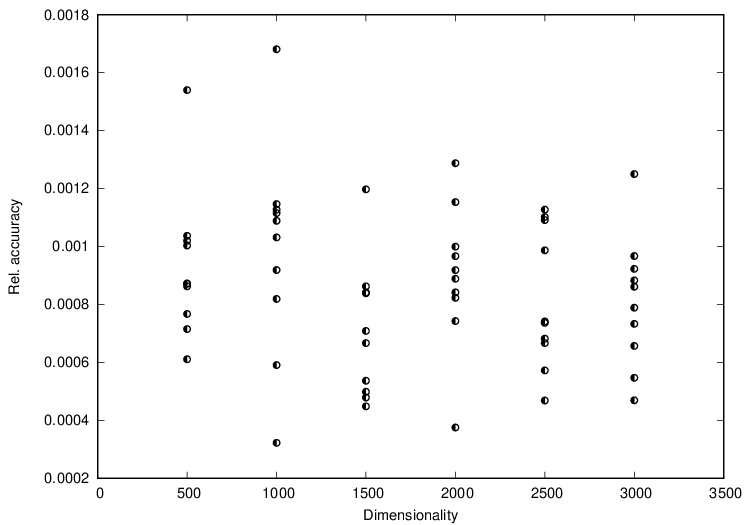}
\caption{Deviations from the optimums}
\label{accu}
\end{figure}

The important characteristic of the quality of solution is its relative accuracy which the projection algorithm manages to attain.
The Fig. \ref{accu} demonstrates the general trend in relative deviation of the objective values obtained by the projection algorithm from optimal values, obtained by \GLPK, and random oscillations in these deviation.
It is worth noticing that we see very little growth in the deviations despite significant growth of the size of problems. Secondly, we see that despite random oscillations the deviations remain quite small, of the order of hundredth of percent.

\begin{figure}
\includegraphics[width=\textwidth]{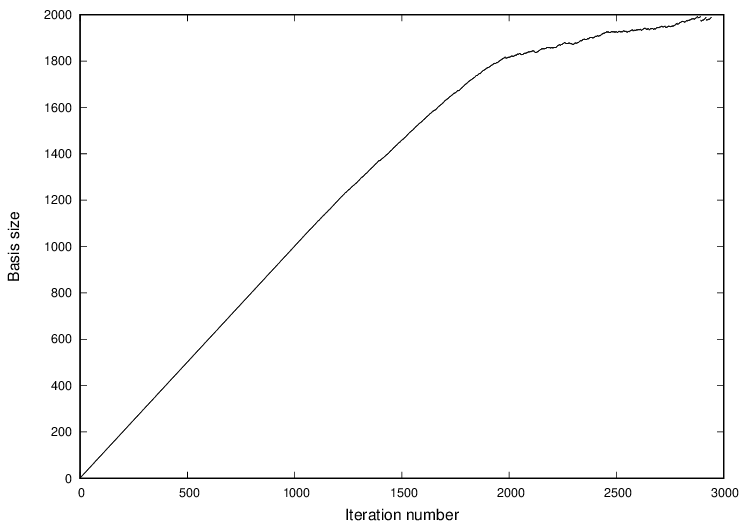}
\caption{Analysis of the projection algorithm: Basis size over iterations.}
\label{copr1}
\end{figure}

\begin{figure}
\includegraphics[scale=0.9]{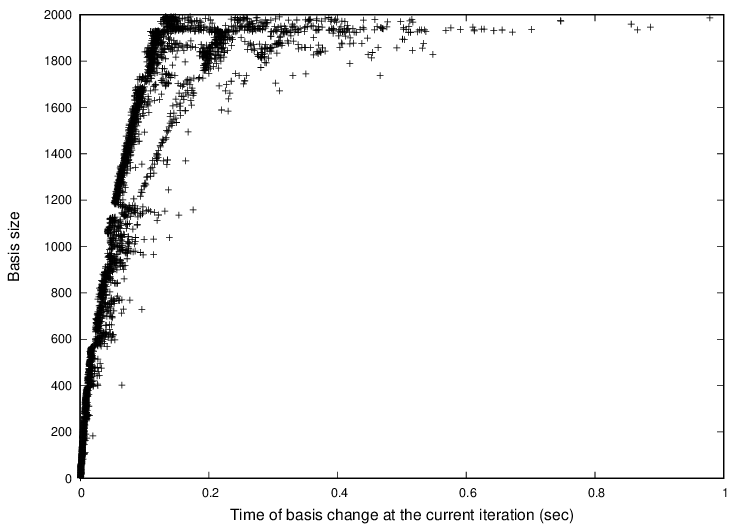}
\caption{Analysis of the projection algorithm: The time it takes to recalculate the basis depending on its size.}
\label{copr2}
\end{figure}

Finally we provide experimental data which give a plausible explanation on why projection algorithms may compete with traditional symplex-like algorithms.
Fig. \ref{copr1} shows the growth of the number of active generators in conical projection. Quite naturally it grows practically linear with algorithm iterations and it is clear that
for the major part of the run the basis size remains well below the theoretical limit. The computing time for update of projection operator in projection algorithm also follows the growth of the basis size (Fig. \ref{copr2}).
It implies that the best part of it essentially smaller then the maximal run time. 

\begin{figure}
\includegraphics[scale=0.9]{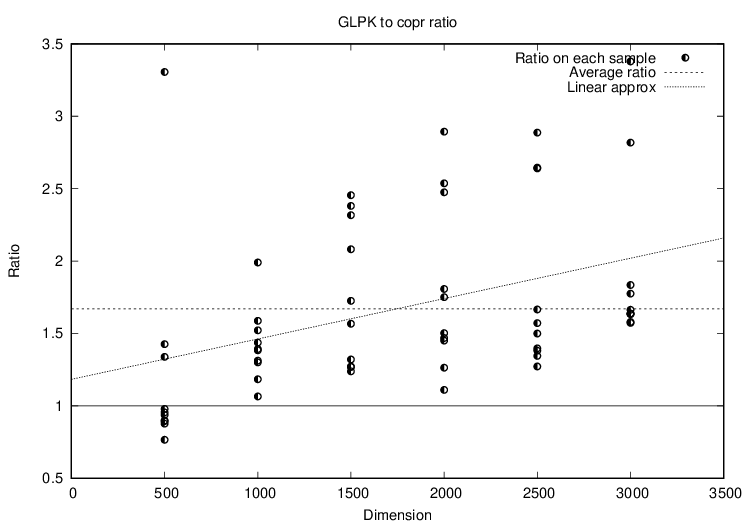}
\caption{The ratio of the running times of \GLPK to copr on each sample.}
\label{ratio}
\end{figure}

Fig. \ref{test_lp1} and \ref{test_lp2} demonstrate generally polynomial growth of computing time both for GLPK and projection algorithm as a function of problems size.
We provide the additional Fig. \ref{ratio} to demonstrate more explicitly the difference in computing time between GLPK and projection algorithm.
It can be seen from Fig. \ref{ratio} that despite rather large fluctuations in the ratio between computing times for GLPK and projection algorithm
the general tendency is in favor of projection when problems become larger.

\newpage

\subsubsection*{Acknowledgements}
This research was initiated and supported at Sirius supplementary educational
program <<Actual methods of information theory and optimization>>, November, 2022.

\end{document}

%% file: packages.tex
\usepackage			{amsfonts}
\usepackage			{amsmath}
\usepackage			{array}
\usepackage			{delarray}
\usepackage[boxed]		{algorithm2e}
\usepackage			{cmap}                   
\usepackage			{comment}
\usepackage			{enumerate}
\usepackage			{graphicx}
\usepackage			{latexsym}
\usepackage			{longtable}
\usepackage[mathscr]		{eucal}
\usepackage			{mathtext}               
\usepackage			{multirow}
\usepackage			{rcs}
\usepackage			{theorem}
\usepackage			{url}
\usepackage			{xspace}

%% file: DeclareMathOperator.tex
\DeclareMathOperator	{\amin}	{{\mathrm argmin}}

\DeclareMathOperator	{\co}	{{\mathrm co}}
\DeclareMathOperator	{\Co}	{{\mathrm Co}}
\DeclareMathOperator	{\Cone}	{{\mathrm Co}}

\DeclareMathOperator	{\lin}	{{\mathrm lin}}

%% file: preamble_ru.tex
\usepackage[T2A]		{fontenc}           	
\usepackage[utf8]		{inputenc}         	
\usepackage[english]	{babel}			

%% file: newcommands.tex
\newcommand{\bleq}[2]{\begin{equation}\label{#1}{#2}\end{equation}}

\newcommand{\R}{\mathbb{R}}

\newcommand{\xs}{{x^\star}}

\newcommand{\cK}{{\cal{K}}}

\newcommand{\half}{\frac12}

\newcommand{\minar}[1]{\begin{array}[t]{c} \min \\ #1 \end{array}}
\newcommand{\maxar}[1]{\begin{array}[t]{c} \max \\ #1 \end{array}}

\newcommand{\vnll}{{\mathbf 0}}
\newcommand{\vone}{{\mathbf 1}}

\newcommand{\prc}[2]{#1\!\downarrow\!#2}
\newcommand{\supp}[2]{{\left({#1}\right)}_{#2}}